\documentclass[10pt]{article}

\usepackage{graphicx}
\usepackage{algorithm}
\usepackage{algorithmic}
\usepackage{subfigure}

\usepackage[a4paper, left=35mm,right=35mm,top=34mm,bottom=34mm]{geometry}
\usepackage[utf8]{inputenc}
\usepackage[T1]{fontenc}
\usepackage[english]{babel}

\usepackage{fancyhdr}

\graphicspath{{/}}

\title{Graph Neural Network-based Surrogate Models for Finite Element Analysis}

\author{Meduri Venkata Shivaditya\thanks{Universit\'e Paris-Saclay, CentraleSup\'elec, MICS, 9 rue Joliot Curie, 91192 Gif-sur-Yvette, France},
Jos\'e Alves\thanks{Transvalor SA, E-Golf Park, 950 Avenue Roumanille, 06410 Biot, France},
Francesca Bugiotti\thanks{Université Paris-Saclay, CNRS, CentraleSup\'elec, LISN, 1 rue Raimond Castaing, 91190 Gif-sur-Yvette, France},
Fr\'ed\'eric Magoul\`es\thanks{Universit\'e Paris-Saclay, CentraleSup\'elec, MICS, 9 rue Joliot Curie, 91192 Gif-sur-Yvette, France, Email: frederic.magoules@hotmail.com}
}

\begin{document}

\maketitle
\thispagestyle{fancy}

\begin{abstract}
Current simulation of metal forging processes use advanced finite element methods.
Such methods consist of solving mathematical equations, which takes a significant amount of time for the simulation to complete.
Computational time can be prohibitive for parametric response surface exploration tasks.
In this paper, we propose as an alternative, a Graph Neural Network-based graph prediction model to act as a surrogate model for parameters search space exploration and which exhibits a time cost reduced by an order of magnitude.
Numerical experiments show that this new model outperforms the Point-Net model and the Dynamic Graph Convolutional Neural Net model. 
\end{abstract}


\section{Introduction}

Finite element methods (FEM) are extensively used to simulate real world physical phenomenon \cite{Cia2002}, \cite{Hug2003}.
In FEM, the mathematical equations associated to the physics are reformulated with a variational formulation which is then discretized.
This discretization is performed on a mesh, and the quality of the elements of the mesh impacts directly the approximated solution.
Several techniques exit to ensure high quality of the approximated solution, including for instance Streamline-Upwind Petrov-Galerkin (SUPG) \cite{Rus2006}, stabilized finite elements \cite{HM2004WM}, bubble elements \cite{DF2008}, infinite elements \cite{Bur1994}, \cite{AM2006}.
Since the FEM solution is obtained by solving a linear system of equations, whose size is proportional to the number of discretization points composing the mesh, this represent a significant part of the computational time.
The high time cost of generating the FEM solution makes it tedious for running thousands of simulations by varying the input parameters for optimization applications and finding the best input parameter set.
This is particularly true in metal forging process design.

In this paper, as an alternative to FEM, we explore deep learning models for metal forging process design.
The motivation behind using a deep learning surrogate model is to create a hybrid approach in which FEM is only used to generate high resolution results in a reduced parametric space.
Neural Networks \cite{csaji2001approximation}
are efficient at learning patterns in data and have been widely used in various applications including image learning,
speech recognition, graph learning and so on, where classical learning is difficult because of the complexity in data.
The low time cost of deep learning model enable faster parameter search space exploration and can save many days' worth of time during the optimization process design.
Here, we created a Graph Neural Network-based deep learning model that takes the mesh objects in the form of a graph and the input parameters as the features of this graph.
The reason for considering Graph Neural Network \cite{refgnn} based approach relies in the property of graphs \cite{xu2018powerful} which share the same permutation in-variance property between meshes and point-cloud objects.
We train the model to some simulations (generated by FEM) to minimize the Mean Squared Error loss for the prediction and the actual graph for all training simulations and test the performance of the model on some test simulations.
The proposed Graph Neural Network (GNN) model, once trained, can generate the output data for a set of inputs in a time smaller than 500 milliseconds, keeping in mind that the actual simulation process using FEM takes 110 minutes to generate the same output data.
The proposed model takes 99.9\% less time than FEM.

This paper is organized as follows,
In section~\ref{section:relwor}, we present some related works for the task of classification and segmentation on point-cloud objects.
In section~\ref{section:method}, we describe in details the methodology.
In section~\ref{sec:building}, we address several issues like framing of the problem, exploring the dataset, evaluation methods, model architectures and results.
In section~\ref{sec:conc}, we present the conclusions.

\section{Related works}
\label{section:relwor}

There are different types of algorithms to create classification and segmentation models, see for instance \cite{ZE2012RSER} and references therein.
Support Vector Machines (SVM) \cite{HDO1998} \cite{ZE2010JACT} \cite{ZE2012JACT} \cite{LML2008IJCM}, Multi-Layer Perceptron \cite{desai2021anatomization}, Convolutional Neural Networks, Recursive Deterministic Perceptron (RDP) \cite{MZE2013EB}, Neural Network, and Deep Learning Methods. 
SVMs construct a hyperplane that separates the objects with the maximum possible margin.
Taking advantage of the Kernel Trick, we are able to project the data into higher dimensions which allows for non-linear margins.
Multi-Layer Perceptrons (MLP), Convolutional Neural Networks (CNN), Deep Learning models all are made up of neurons.
MLPs are stacked fully connected layers where a neuron of a layer is connected to every neuron in the other layer.
In CNNs, instead of matrix multiplications, convolutional operation is employed.
The Point-Net model \cite{QCH2017} is one of the first iterations of deep learning models which can deal with permutation invariant data structures like the point cloud object.
The Point-Net model consists of various 1D convolutional layers and activation functions followed by a Global Max Pooling with fully connected layers in the end, modified for the objective of the model. 
The Dynamic-Graph Convolutional Neural Network model forms connections between the points of the point-cloud object after every layer following a $K$-Nearest Neighbors algorithm.
This graph structure allows for pooling of local level features and takes the neighborhood into consideration.
This approach drastically improved the performance of Segmentation tasks on Point-Cloud data.
Graph Neural Network-based models were applied on 3D objects for the task of object classification and segmentation using Dense Graph Convolutional layers.
These approaches were mainly developed for classification/segmentation and not for data prediction purposes. 

FEM simulate physical phenomena modeled by Partial Differential Equations (PDE) by solving a discrete problem \cite{Cia2002}, \cite{MZ2018CMAME}, \cite{AM2007}.
There have been attempts to approximate the PDE solutions using deep learning \cite{belbute2020combining}.
These models have been tested on complex fluid dynamic simulations and use a hybrid (graph) neural network that combines a traditional graph convolutional network with an embedded differentiable fluid dynamics simulator inside the network itself for faster results and a good approximation. 

The 3D objects considered is usually in the point-cloud format which ignores the topology of the object.
The Point-Net, DGCNN and GNN models were not commonly used as surrogates for simulation procedures but were used for automatic analysis of data generated by FEM.
In the following section, we explore deep learning models using Graph Neural Networks for data prediction.

\section{Methodology}
\label{section:method}

\subsection{Dataset}

We used data arising from FEM simulations of a Yoke metal forging process. 
Generating each simulation on a quad-core machine takes approximately 110 minutes, which appears to be a very expensive time cost for parameter search space exploration and optimization tasks. 
The input parameters provided to the simulation are initial temperature of the billet and friction coefficient between the billet and the lower deformable die. 
The goal is to predict the wear prediction as a surface field on the lower deformable die/upper deformable die for a given set of initial input conditions.
The proposed deep learning graph neural network model will approximate the final lower deformable die/upper deformable die for the given input set.
The input parameters of temperature and friction coefficient are added as node features to each of the node of the mesh.
The surface Lower Deformable Die (LDD) consists of 9 215 nodes.
The surface Upper deformable die (UDD) consists of 6 617 nodes.
Each Cell of the mesh has features like the normal stress, flow stress and wear, but all the features that cannot change have been eliminated to reduce the complexity of the model. 
The points of the mesh and their features are considered as the node features of the graph.
Within the total number of 40 simulations in the dataset, 30 simulations are considered for training the proposed model and 10 simulations are considered for testing the proposed model.

\subsection{Finite element mesh to graph conversion}

\begin{figure}[h]
    \centering
    \includegraphics[width=0.28345\linewidth]{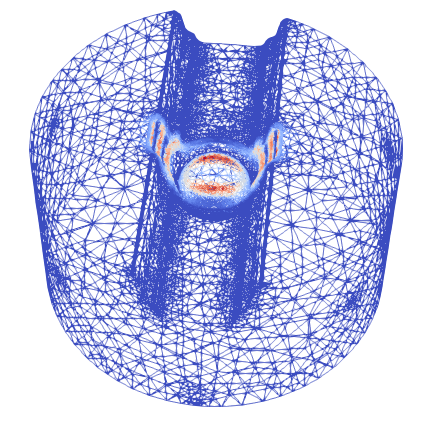}	
    \includegraphics[width=0.28345\linewidth]{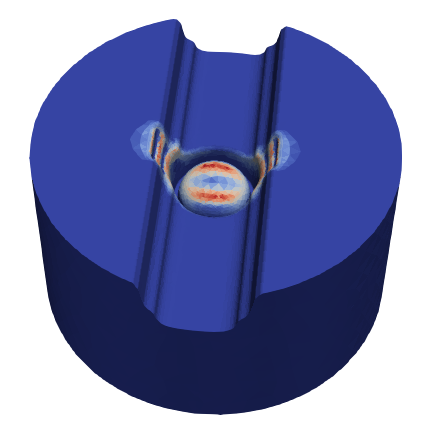}
    \caption{Lower Deformable Die. Mesh (left) and finite element solution (right).}
		\label{fig:vtk}
\end{figure}

The model should take as input the starting mesh object and the set of input parameters (temperature, friction coefficient, \ldots) and as output the final simulation object with wear field for each of the nodes. 
The mesh object is converted to graph and the input parameters are added to every node of the graph as node features.
We extract the mesh topological connectivity and use it to connect nodes in the corresponding graph data structure.
The connections are between the points and the features are associated with the cells of the mesh.
So, we convert the cell features to point features by averaging the feature values of all the cells the point belongs to.
An example of the visualization \cite{MP2007} is shown in Figure~\ref{fig:vtk}.

\subsection{Mathematical formulation}

We design a deep learning model that directly takes a graph as an input and returns a graph as output.
A graph is represented as a set of node features, edge features, graph features and node-to-node connectivity.
For each simulation, we will convert the mesh to a graph by taking the point features as the node features, and the point-to-point connectivity of the mesh as the node-to-node connectivity of the graph.
The proposed model will predict the final simulation in the form of a graph with node features as the point features.
The proposed model once trained can be used as a surrogate model for FEM.
Since we have independent graphs for training, testing and validation, we will use the inductive setting \cite{NIPS2017_5dd9db5e} of graph models for this problem.
Because we are using the inductive setting, we don't have to worry about cross talk or information leakage from unseen data.

\section{Building the model}
\label{sec:building}

The required steps to build the model include: (i) processing the Dataset for training and testing, (ii) defining the evaluation metric to test the model's performance, (iii) defining the model's architecture, (iv) training and testing the model's performance with multiple datasets and other relevant models.
Below we describe all the steps to build the GNN Surrogate model we proposed, which will predict what the final result will be, and we evaluate it's performance.

\subsection{Evaluation metric}

Since, the task is a graph node regression problem, we will use the Mean Squared Error (MSE) as a metric to compute the loss of our model.
For each node of the output graph, we calculate the difference squared with the actual output and average it over all the nodes.
This will be the MSE loss for the first simulation, which will be used to readjust the weights of the model.  
We trained the model to minimize the Mean Squared Error loss.
After completely training the models, we tested our model on the unseen data using the MSE loss metric.  
$MSE = \frac{1}{n}\sum_{i=1}^{n}(\hat{y}_i-y_i)^2$ where $n$ is the number of nodes in the graph, $y_i$ is the actual wear for node $i$ in the graph, and $\hat{y}_i$ is the predicted wear for node $i$ in the graph.
We will also use the Root Mean Squared Error (RMSE) and the coefficient of determination to evaluate how well the trained model is performing on the test data.
Coefficient of Determination $R^2$ is the proportion of variation in the target variable $y_i$ that is explained by the model and is calculated as follows: $R^2 = 1 - \frac{SS_{res}}{SS_{tot}}$, where $SS_{res} = \sum_{i=1}^{n} \epsilon_i^2$, $SS_{tot} = \sum_{i=1}^{n} (y_i - \hat{y})^2$ with $\epsilon_i$ equal to $\hat{y}_i-y_i$.
RMSE and $R^2$ are much more comprehensible to the human eye when assessing the performance of the model.

\subsection{Graph Neural Network-based surrogate model for FEM simulations}

Neighbourhood information of nodes in graphs plays an important role in prediction of the target feature of the nodes.
The problem of considering neighbourhood information through convolutions to make predictions for graphs is unlike a structured object like images where the pixels are arranged in a certain structure.
Graph Convolutional Layers are excellent approaches to take the neighbourhood information to make predictions using deep learning models for graph like permutation invariant data structures. 
The proposed model consists of five Graph Convolutional Layers followed by Re-LU activation layers.
To avoid over-fitting, Dropout method of regularization is introduced with a probability of dropout of 1\%.
The model is designed to generate the wear field, and since wear can only take positive values, the proposed model is designed with a Re-LU layer as the final layer which ensures that the model outputs are with positive values for each of the nodes.

\subsection{Graph convolutional layer}

Each node of the graph constitutes a computation graph built using it's neighbours \cite{zhang2019graph}.
The neighbourhood information is aggregated using a permutation invariant operation like average or maximum, after which a neural network is applied on the aggregated information before sending it to the next layer.
This process of information propagation through the neighbourhood to the central node ensures that the information passed is not altered if the permutations are changed for the same computation graph.
Each graph convolutional layer takes in node features for information propagation, and the neural networks of the layers control the number of node features for the output graph.
For the first graph convolutional layer in the proposed GNN surrogate model, the layer takes a graph with 5 node features as input and returns a graph with 50 node features as output before applying the activation function.
The activation functions add non-linearity to the model.
The second graph convolutional layer takes a graph with 50 node features and returns a graph with 100 node features followed by Re-LU activation layer.
Then, in the following 3 graph convolutional layers with activation functions, the 100 features are reduced to 50 features and then to 1 feature which is the target wear feature to which the Mean Squared Error loss function is applied.
The dropout layers are used for regularization and are stacked after each activation function.
The final layer is a Re-LU layer to ensure that the model only predicts positive values for the nodes.

\subsection{Training the model}

\begin{figure}[h]
    \centering
    \includegraphics[width=0.3545\linewidth]{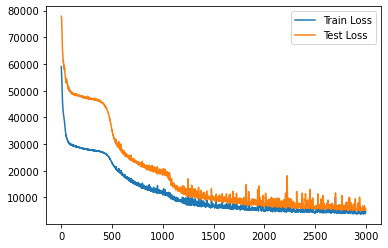}
    \includegraphics[width=0.3545\linewidth]{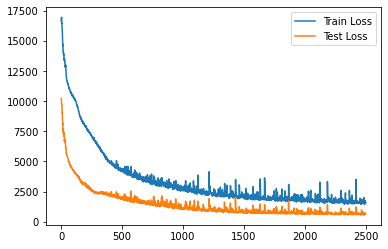}
    \caption{Left: Loss Curve (LDD). Right: Loss Curve (UDD)}
		\label{fig:loss}
\end{figure}

The model's architecture is trained using the deep graph library framework \cite{wang2019dgl} based on Pytorch. 
We trained the models on two separate datasets predicting the final mesh of the lower deformable die and the upper deformable die.
The models are trained for 2 500 epochs with
a method for stochastic optimization using a learning rate of $8\times10^{-4}$.
Each of the datasets consists of 40 simulations, where 30 of which are used for training and 10 are used for testing.
The Data-Loader is defined with a batch size of 1, i.e., only 1 simulation is used per iteration of model training. 
The loss-curves, illustrated in Figure~\ref{fig:loss}, suggest that the models trained on LDD and UDD datasets are properly trained without any over-fitting.
After 2 500 epochs, the train and test loss per simulation is nearly the same, suggesting that the model is well trained.

\subsection{Bench-marking the new model}

The model's prediction seems to be a good reflection of the wear.
After training, the model obtained 73.56 N/m RMSE loss on average per simulation on the test LDD dataset. On the test UDD dataset, the model obtained 26.44 N/m RMSE loss on average. Average Coefficient of Determination on the test data for LDD is 93.7\% and for UDD is 92.8\%.

Along with the proposed model, we also tested Point-Net and Dynamic Graph Convolutional Neural Network models using the 75\%/25\% split model experimentation methodology.
Point-Net segmentation model's last layers have been modified to predict a continuous value for each of the nodes. 
Dynamic Graph Convolutional Neural Net segmentation model's last layers are also modified to accommodate for the dataset.
For the DGCNN model, we chose dynamic neighbours as 5 and the embedding dimension size as 64. 
Early Stopping mechanism is implemented while training to avoid over-fitting.

The Point-Net and the DGCNN models take point-cloud objects as input and not the node-to-node connectivity. 
Point-Net model performs 1d convolution operations and activations for several layers before performing a Global Max Pooling to aggregate the features, which are later fit into a fully connected layer to output the prediction. 
There is no necessity for node-to-node edges for the Point-Net model. The DGCNN model computes the connections between nodes using the $K$-Nearest Neighbours method and updates connections after each layer. 
The proposed GNN model has the advantage of taking the actual connections of the topological structure of the mesh and aggregate neighbourhood information using that.

\begin{table}[h]
\centering
  \caption{RMSE for the models tested}
  \label{tab:freq}
	\small
  \begin{tabular}{lrr}
    Model&LDD dataset&UDD Dataset\\
		\hline
		\hline
    Point-Net model & 299.11 & 137.58\\
    DGCNN model & 303.31 & 126.79\\
    Proposed GNN model & {73.56} & {26.43}
\end{tabular}
\end{table}

The values shown in Table~\ref{tab:freq} correspond to the root mean squared error loss per node of the final data predicted on average on the testing data.
All the models are tested on two datasets corresponding to the Lower Deformable Die and the Upper Deformable Die respectively.
This is used to compare the models' performance.

The GNN model performs extremely well when compared to the other models.
The RMSE loss on the testing dataset is significantly better than the other models.
The point-cloud based deep learning models didn't fare well in comparison with the GNN model.
The Point-Net model works well when compared to DGCNN, but couldn't compete with the GNN model.
The GNN's ability to take the neighbourhood information using the topological connections of the mesh might be the reason why the proposed GNN model outperforms the Point-Net model and the DGCNN model.

\subsection{Time cost analysis using the proposed surrogate model}

FEM simulations were generated on a quad-core computer with 16 CPUs.
A single simulation took about 1 hour 50 minutes (on 4 CPUs/case), but the 40 cases were automatically launched on 16 CPUs meaning that 4 cases were running simultaneously in parallel.
The total effective time for generating 40 simulations is 18 hours.
After training the surrogate model on the same CPUs,
the average time taken to generate the final result is 500 milliseconds, this means that using the GNN Surrogate, we can save the time cost by 99.9\%.

\section{Conclusions}
\label{sec:conc}

Using the Graph Neural Network approach to act as surrogate model for FEM is a well performing approach.
Depending, on the time-cost versus accuracy required, we can adapt the model ensuring flexibility.
An important point to mention is that the surrogate model we proposed performs well, even though it was trained only on 30 simulations data.
The Graph Neural Network surrogate model, because of its cheap time cost, can aid in parameter search space exploration.

\section{Acknowledgements}
The authors would like to thank M.~Yu, R.~Raj, W.~He, Z.~Liu and Y.~Chen for the useful discussions.

\end{document}